\documentclass{amsart}

\usepackage{amsfonts}
\usepackage{amsmath}
\usepackage{amsthm}
\usepackage{amssymb}
\usepackage[english]{babel}
\usepackage{graphics}
\usepackage{latexsym}
\usepackage{longtable} 
\usepackage{mathrsfs} 
\usepackage{graphicx}
\usepackage{wrapfig} 
\usepackage[pdftex,colorlinks=true]{hyperref}

\newtheorem{theorem}{Theorem}
\newtheorem{corollary}[theorem]{Corollary}
\newtheorem{lemma}[theorem]{Lemma}
\newtheorem{proposition}[theorem]{Proposition}

\newtheorem{claim}[theorem]{Claim}
\newtheorem{example}[theorem]{Example}

\theoremstyle{definition}
\newtheorem{definition}[theorem]{Definition}
\newtheorem{remark}[theorem]{Remark}

\newcommand{\mK}{\mathbb{K}}
\newcommand{\mW}{\mathbb{W}}

\newcommand{\mP}{\mathscr{P}}
\newcommand{\mM}{\mathcal{M}}

\newcommand{\mY}{\mathscr{Y}}
\newcommand{\mD}{\mathcal{D}}

\newcommand{\R}{\mathbb{R}}

\newcommand{\N}{\mathbb{N}}

\newcommand{\mB}{\mathbb{B}}
\newcommand{\X}{\textbf{X}}

\newcommand{\um}{{\underline{m}}}
\newcommand{\uX}{\underline{\textbf{X}}}

\newcommand{\noi}{\noindent}
\newcommand{\ms}{\medskip}

\newcommand{\al}{\alpha}
\newcommand{\be}{\beta}
\newcommand{\ga}{\gamma}

\newcommand{\de}{\delta}
\newcommand{\De}{\Delta}
\newcommand{\e}{\varepsilon}

\newcommand{\Om}{\Omega}
\newcommand{\om}{\omega}

\newcommand{\lharpoonup}{-\!\!\!\!-\!\!\!\!\rightharpoonup}

\newcommand{\larrow}{\longrightarrow}
\newcommand{\ot}{\otimes}

\newcommand{\lmapsto}{\longmapsto}
\newcommand{\ri}{\rightarrow}

\newcommand{\p}{\partial}
\newcommand{\sub}{\subseteq}
\newcommand{\set}{\setminus}
\newcommand{\by}{\times}

\newcommand{\ess}{\textrm{ess}}

\newcommand{\spn}{\textrm{span}}
\newcommand{\supp}{\textrm{supp}}

\newcommand{\bt}{\begin{theorem}}\newcommand{\et}{\end{theorem}}
\newcommand{\bd}{\begin{definition}}\newcommand{\ed}{\end{definition}}
\newcommand{\bl}{\begin{lemma}}\newcommand{\el}{\end{lemma}}
\newcommand{\beq}{\begin{equation}}\newcommand{\eeq}{\end{equation}}
\newcommand{\bc}{\begin{claim}}\newcommand{\ec}{\end{claim}}
\newcommand{\bex}{\begin{example}}\newcommand{\eex}{\end{example}}
\newcommand{\bcor}{\begin{corollary}}\newcommand{\ecor}{\end{corollary}}
\newcommand{\bp}{\begin{proof}}\newcommand{\ep}{\end{proof}}

\newcommand{\BPL}{\medskip \noindent \textbf{Proof of Lemma} }

\newcommand{\BPCOR}{\medskip \noindent \textbf{Proof of Corollary} }

\newcommand{\BPT}{\medskip \noindent \textbf{Proof of Theorem} }

\numberwithin{equation}{section}

\begin{document}

\title[Mollification of $\mathcal{D}$-solutions]{Mollification of $\mathcal{D}$-solutions to Fully Nonlinear PDE Systems}

\author{Nikos Katzourakis}
\address{Department of Mathematics and Statistics, University of Reading, Whiteknights, PO Box 220, Reading RG6 6AX, Berkshire, UK}
\email{n.katzourakis@reading.ac.uk}



\date{}



\begin{abstract} In a recent paper the author has introduced a new theory of generalised solutions which applies to fully nonlinear PDE systems of any order and allows the interpretation of merely measurable maps as solutions. This approach is duality-free and builds on the probabilistic representation of limits of difference quotients via Young measures over certain compactifications of the ``state space". Herein we establish a systematic regularisation scheme of this notion of solution which, by analogy, is the counterpart of the usual mollification by convolution of weak solutions and of the mollification by sup/inf convolutions of viscosity solutions.
\end{abstract}

\maketitle


\section{Introduction} \label{section1}

Let $p,n,N,M \in \N$, $\Om\sub \R^n$ an open set and 
\beq \label{1.1}
F\ : \ \Om \by \Big(\R^N\by \R^{Nn} \by \cdots \by \R^{Nn^p}_s \Big) \larrow \R^M
\eeq
a Carath\'eodory mapping. Here $\R^{Nn}$ denotes the space of $N\!\by\! n$ matrices wherein the gradient matrix 
\[
Du(x)\ =\ \big( D_iu_\al(x)\big)_{i=1,...,n}^{\al=1,...,N}
\]
of (smooth) maps $u:\Om \sub \R^n \larrow \R^N$ is valued, whilst $\R^{Nn^p}_s$ denotes the space of symmetric tensors 
\[
\begin{split}
&\Big\{ \X \in \R^{Nn^p}\ \big| \ \X_{\al i_1 ... i_a...i_b... i_p} = \X_{\al i_1 ... i_b...i_a... i_p}\ , \  \\
&\ \, \al =1,...,N, \ i_1,...,i_p \in\{1,...,n\},\ 1\leq a\leq b\leq p\Big\}
 \end{split}
\]
\
wherein the $p$th order derivative 
\[
D^p u(x)\ =\ \left( D^p_{i_1 ... i_p}u_\al(x)\right)_{i_1,...,i_p \in \{1,...,n\}}^{\al=1,...,N}
\]
is valued. Obviously, $D_i\equiv \p/\p x_i$, $x=(x_1,...,x_n)^\top$, $u=(u_1,...,u_N)^\top$ and $\R^{Nn^1}_s= \R^{Nn}$. In the recent paper \cite{K8} we introduced a new theory of generalised solutions which allows for merely measurable maps to be rigorously interpreted and studied as solutions of systems with even discontinuous coefficients and without requiring any structural assumptions (like ellipticity or hyperbolicity). Namely, our approach applies to measurable solutions  $u:\Om \sub \R^n \larrow \R^N$ of the $p$th order system
\beq \label{1.2}
F\Big(x,u(x),D^{[p]}u(x)\Big)\, =\, 0,\ \  \ x\in \Om,
\eeq
where 
\[
D^{[p]} u\, :=\, \big(Du,D^2u,...,D^pu\big)
\]
denotes the \emph{$p$th order Jet of $u$}. Since we do not assume that solutions must be locally integrable on $\Om$, the derivatives $Du, ...,D^pu$ may not have a classical meaning, not even in the distributional sense. Using this new approach, in the very recent papers \cite{K8,K9,K10} we studied efficiently certain interesting problems arising in PDE theory and in vectorial Calculus of Variations which we discuss briefly at the end of the introduction. 

In the present paper we are concerned with the development of a systematic method of mollification of generalised solutions to the fully nonlinear system \eqref{1.2} by constructing approximate smooth solutions to approximate systems. The mollification method we establish herein is the counterpart of the standard mollification by convolution which is the standard analytical tool in the study of weak solutions and to the so called sup/inf convolutions used in the theory of viscosity solutions of Crandall-Ishii-Lions (for a pedagogical introduction we refer to \cite{K7}).

Our starting point for the definition of solution is not based either on standard duality considerations via integration-by-parts (the cornerstone of weak solutions) or on the maximum principle (the mechanism of the more recent method of viscosity solutions). Instead, we build on the probabilistic representation of limits of difference quotients by utilising \emph{Young} measures, also known as parameterised measures. These are well-developed objects of abstract measure theory of great importance in Calculus of Variations and PDE theory (see e.g.\ \cite{K8} and \cite{E, P, FL, CFV, FG, V, KR}). In the present setting, a version of Young measures is utilised in order to define generalised solutions of \eqref{1.2} by applying them to the difference quotients of the candidate solution. The essential idea restricted to the first order case $p=1$ of \eqref{1.2} goes as follows: let $u\in W^{1,1}_{\text{loc}}(\Om,\R^N)$ be a strong solution to 
\beq   \label{1.3}
F\big(x,u(x),Du(x) \big)\,=\, 0, \quad \text{ a.e.\ }x\in \Om.
\eeq
Let us now rewrite \eqref{1.3} in the following unconventional fashion:
\[
\sup_{X\in \, \supp(\de_{D u(x)})} \big| F\big( x,u(x),X\big)\big|\, =\, 0, \quad \text{ a.e. }x\in \Om.
\]
That is, we understand the gradient $Du : \Om\sub \R^n  \larrow\R^{Nn}$ as a probability-valued map given by the Dirac measure at the gradient 
\[
\de_{Du} \ :\ \Om \sub \R^n  \larrow \mP(\R^{Nn}), \ \ \ x\lmapsto \de_{Du(x)}, 
\]
in hopes of relaxing the requirement to have a concentration measure. The goal is to allow instead general probability-valued maps arising as limits of the difference quotients for nonsmooth maps. Indeed, if $u : \Om \sub \R^n \larrow \R^N$ is only \emph{measurable}, we consider the probability-valued mappings 
\[
\de_{D^{1,h}u} \ : \ \Om \sub \R^n \larrow \mathscr{P}\big( \smash{\overline{\R}}^{Nn} \big) , \ \ \ x\lmapsto \de_{D^{1,h}u(x)},
\]
where $D^{1,h}$ is the usual difference quotients operator and $ \smash{\overline{\R}}^{Nn}$ is the Alexandroff 1-point compactification $\R^{Nn}\cup\{\infty\}$. Namely, we view $D^{1,h}u$ as an element of the space of Young measures $\mY\big( \Om, \smash{\overline{\R}}^{Nn} \big)$ (see the next Section \ref{section2} for the precise definitions). By using that $\mY\big( \Om, \smash{\overline{\R}}^{Nn} \big)$ is weakly* compact, there always exist probability-valued mappings $\mD u \in \mY\big( \Om, \smash{\overline{\R}}^{Nn} \big)$
such that along infinitesimal subsequences $(h_\nu)_1^\infty$ we have
\beq  \label{1.4}
\de_{D^{1,h_\nu}u} \, \overset{^*}{\smash{\lharpoonup}}\, \, \mD u \quad \text{ in }\, \mY\big( \Om, \smash{\overline{\R}}^{Nn} \big), \ \ \ \text{ as }\nu \ri \infty
\eeq
(even if $u$ is merely measurable). Then, we require\footnote{The version of the definition we are using herein is different from the one we put foremost in \cite{K8}-\cite{K10} because this simplifies the proofs that follow. In \cite{K8} we proved several equivalent formulations of the same notion which we will not utilise, so we take this version as primary here.}
\beq  \label{1.5}
\sup_{X\in \, \supp_*(\mD u(x))} \big| F\big( x,u(x),X\big)\big|\, =\, 0, \quad \text{ a.e. }x\in \Om,
\eeq
 for any``diffuse" gradient $\mD u$, where 
 \[
 \supp_*(\mD u(x))\, := \, \supp(\mD u(x))\set\{\infty\}.
 \]
 Since \eqref{1.4} and \eqref{1.5} are \emph{independent} of the regularity of $u$, they can be taken as a notion of \textbf{diffuse derivatives} of the measurable map $u : \Om \sub \R^n \larrow \R^N$ and of \textbf{$\mD$-solutions} to the PDE system \eqref{1.3} respectively. If $u$ happens to be weakly differentiable, then we have $\mD u=\de_{Du}$ a.e.\ on $\Om$ and we reduce to strong solutions. Except for a small further technical generalisation (we may need to take special difference quotients depending on $F$), \eqref{1.4} and \eqref{1.5} comprise our notion of generalised solutions in the first order case of \eqref{1.2}.
 
This paper is organised as follows. The introduction is followed by Section \ref{section2} which is a quick review of the main points of \cite{K8} necessary for this work. The main results of this paper are in Section \ref{section3} (Theorem \ref{theorem11} and Corollary \ref{corollary12}) and establish the main properties of our approximations. The technical core of our analysis is contained in the preparatory Lemma \ref{lemma10}. We expect the analytical results established herein to play a prominent role in future developments of the theory, but we refrain from providing any immediate applications in this paper.
 
We conclude this introduction with some results recently obtained by using the technology of $\mD$-solutions. Our motivation to introduce them primarily comes from the necessity to study the recently discovered equations arising in vectorial Calculus of Variations in the space $L^\infty$, that is for variational problem related to functionals 
\beq  \label{1.6}
E_\infty(u,\Om)\, :=\, \big\| H(\cdot,u,Du) \big\|_{L^\infty(\Om)}
\eeq
applied to Lipschitz maps $u : \Om \sub \R^n \larrow \R^N$ (for an introduction to the topic we refer to \cite{C,BEJ,K7}). In the simplest case of $H(\cdot,u,Du)=|Du|^2$, the analogue of the Euler-Lagrange equation is the $\infty$-Laplace system:
\beq \label{1.7}
\De_\infty u \, :=\, \Big(Du \ot Du + |Du|^2[Du]^\bot \! \ot I \Big):D^2u\, =\, 0,
\eeq
where  $[Du]^\bot:=\text{Proj}_{(R(Du))^\bot}$. In index form \eqref{1.7} reads
\[
\begin{split}
\sum_{\be=1}^N\sum_{i,j=1}^n \Big(D_iu_\al \, D_ju_\be +\, |Du|^2 & [Du]_{\al \be}^\bot \, \de_{ij}\Big)\, D_{ij}^2u_\be\, =\, 0, \ \ \ \ \al=1,...,N
\end{split}
\]
and $[Du]^\bot$ is the orthogonal projection on the orthogonal complement of the range of $Du$. The vectorial case of the theme has been pioneered by the author in a series of recent papers \cite{K1}-\cite{K6}, while the scalar case is relatively standard by now and has been pioneered by Aronsson in the 1960s (\cite{K7}). In the paper \cite{K8} we studied the Dirichlet problem for \eqref{1.7}, while in \cite{K9} we studied the Dirichlet problem for the system arising from the general functional \eqref{1.6} for $n=1$. In \cite{K9} we also considered the Dirichlet problem for fully nonlinear 2nd order degenerate elliptic systems and in \cite{K10} we considered the problem of equivalence between distributional and $\mD$-solutions for linear symmetric hyperbolic systems.

\section{A quick guide to $\mD$-solutions for fully nonlinear systems} \label{section2}

\subsection{Preliminaries}  \label{subsection2.1} We begin with some basic material needed in the rest of the paper. 

\ms

\noi \textbf{Basics.} The constants $n,N\in \N$ will always denote the dimensions of the domain and the target of our candidate solutions $u:\Om\sub \R^n \larrow \R^N$ defined over an open set. Such mappings will always be understood as being extended by zero on $\R^n\set \Om$. Unless indicated otherwise, Greek indices $\al,\be,\ga,...$ will run in $\{1,...,N\}$ and latin indices $i,j,k,...$ (perhaps indexed $i_1,i_2,...$) will run in $\{1,...,n\}$, even when their range of summation may not be given explicitly. The norm symbols $|\cdot|$ will always mean the Euclidean ones, whilst Euclidean inner products will be denoted by either ``$\cdot$" on $\R^n,\R^N$ or by ``$:$" on tensor spaces. For example, on $\R^{Nn^p}_s$ we have 
\[
|\X|^2 \,= \, \sum_{\al,i_1,...i_p}\X_{\al i_1...i_p}\X_{\al i_1...i_p}\, \equiv \, \X :\X.
\]
Our measure theoretic and function space notation is either standard as e.g.\ in \cite{E,E2} or self-explanatory. For example, the modifier ``measurable" will always mean ``Lebesgue measurable", the Lebesgue measure will be denoted by $|\cdot|$, the characteristic function of a set $E$ by $\chi_E$, the $L^p$ spaces of maps  $u : \Om \sub \R^n \larrow  \R^N$ by $L^p(\Om,\R^N)$, etc.  

We will systematically use the Alexandroff $1$-point compactification of $\R^{Nn^p}_s$. The metric topology on it will be the standard one which makes it isometric to the sphere of equal dimension (via the stereographic projection which identifies the the north pole with infinity $\{\infty\}$). It will denoted by 
\[
\smash{\overline{\R}}^{Nn^p}_s\, :=\ \R^{Nn^p}_s \cup \{\infty\}.
\]
We also note that balls taken in $\R^{Nn^p}_s$ (which we will view as a metric vector space isometrically contained into $\smash{\overline{\R}}^{Nn^p}_s$) will be understood as the Euclidean. 

\ms

\noi \textbf{Young Measures.} Let $E$ be a measurable subset of $\R^n$ and $\mK$ a compact subset of some Euclidean space, which we will later take to be $\smash{\overline{\R}}^{Nn}\!\by \cdots\by \smash{\overline{\R}}^{Nn^p}_s$.

\begin{definition}[Young Measures] The set of Young Measures $\mY(E,\mK)$ consists of the  probability-valued mappings
\[
\vartheta\ : \  E \larrow \mP(\mK), \ \ \ x\lmapsto \vartheta(x),
\]
which are measurable in the following weak* (i.e.\ pointwise) sense: for any continuous function $\Psi \in C^0(\mK)$, the function $E\sub \R^n \larrow \R$ given by
\[
\int_{\mK} \Psi(X)\, d[\vartheta(\cdot)](X)\ : \ \ \ x\lmapsto \int_{\mK} \Psi(X)\, d[\vartheta(x)](X) 
\]
is (Lebesgue) measurable.
\end{definition}
The set $\mY(E,\mK)$ can be identified with a subset of the unit sphere of a certain $L^\infty$ space and this provides very useful compactness and other properties. Consider the $L^1$ space of Bochner integrable maps
\[
L^1\big( E, C^0(\mK)\big)
\]
which are valued in the separable space $C^0(\mK)$ of continuous functions over $\mK$. For background material on these spaces we refer e.g.\ to \cite{FL, Ed, F, V}. The elements of $L^1\big( E, C^0(\mK)\big)$ coincide with the Carath\'eodory functions
\[
\Phi\ :\ E \by \mK \larrow \R, \quad (x,X)\mapsto \Phi(x,X)
\]
which satisfy
\[
\| \Phi \|_{L^1( E, C^0(\mK))}\, :=\, \int_E  \big\|\Phi(x,\cdot)\big\|_{C^0(\mK)}\, dx\, <\, \infty
\]
in the sense that each such $\Phi$ induces a map $E\ni x\mapsto \Phi(x,\cdot) \in C^0(\mK)$. By Carath\'eodory functions we mean that for every $X\in \mK$ the function $x\mapsto\Phi(x,X)$ is measurable and for a.e.\ $x\in E$ the function $X\mapsto \Phi(x,X)$ is continuous. The Banach space $L^1\big( E, C^0(\mK)\big)$ is separable and by using the duality $\big(C^0(\mK)\big)^*=\mM(\mK)$, it can be shown that (see e.g.\ \cite{FL})
\[
\left( L^1\big( E, C^0(\mK)\big) \right)^* \, =\, L^\infty_{w^*}\big( E,\mM(\mK) \big).
\]
The dual space $L^\infty_{w^*}\big( E,\mM(\mK) \big)$ consists of measure-valued maps
\[
E\ \ni \ x \longmapsto \vartheta(x)  \ \in\, \mM(\mK)
\]
which are weakly* measurable and the norm of the space is given by
\[
\| \vartheta \|_{L^\infty_{w^*} ( E,\mM(\mK) )}\, :=\, \underset{x\in E}{\ess\,\sup}\, \left\|\vartheta(x) \right\|(\mK).
\]
Here ``$\|\cdot\|(\mK)$" denotes the total variation on $\mK$. The duality pairing between the spaces
\[
\langle\cdot,\cdot\rangle\ :\  \ L^\infty_{w^*}\big( E,\mM(\mK) \big) \by L^1\big( E, C^0(\mK)\big) \larrow \R
\]
is given by
\[
\langle \vartheta, \Phi \rangle\, :=\, \int_E \int_{\mK} \Phi(x,X)\, d[\vartheta(x)] (X)\, dx.
\]
Then, the set of Young measures can be identified with a subset of the unit sphere of $L^\infty_{w^*}\big( E,\mM(\mK) \big)$:
\[
\mY(E,\mK)\, =\, \Big\{ \vartheta\, \in \, L^\infty_{w^*}\big( E,\mM(\mK) \big)\, : \, \vartheta(x) \in \mP(\mK),\text{ for a.e. }x\in E\Big\}.
\]
Since $L^1\big( E, C^0(\mK)\big)$ is separable, the unit ball of $L^\infty_{w^*}\big( E,\mM(\mK) \big)$ is sequentially weakly* compact. Hence, for any bounded sequence $(\vartheta^m)_1^\infty \sub L^\infty_{w^*}\big( E,\mM(\mK) \big)$, there is a limit map $\vartheta$ and a subsequence of $m$'s along which $\vartheta^m \overset{^*}{\smash{\lharpoonup}}\, \vartheta$ as $m\ri \infty$. 

\begin{remark}[Properties of $\mY(E,\mK)$] \label{remark2} The following facts about Young measures will be extensively used hereafter (the proofs can be found e.g.\ in \cite{FG}): 
\smallskip

i) [\textbf{Weak* compactness of Y.M.}] The set of Young measures is convex and sequentially compact in the weak* topology induced from $L^\infty_{w^*}\big( E,\mM(\mK) \big)$.  

\smallskip
ii)  [\textbf{Mappings as Y.M.}]  The set of measurable maps $v : E\sub \R^n \larrow \mK$ can be imbedded into $\mY(E,\mK)$ via the mapping $v \mapsto \de_v$ given by $\de_v(x):= \de_{v(x)}$.  
\end{remark}

The next lemma is a minor variant of a classical result (see \cite{K8, FG, FL}) but it plays a fundamental role in our setting because it guarantees the compatibility of strong solutions with $\mD$-solutions. 

\begin{lemma} \label{lemma2} Suppose $E\sub \R^n $ is a measurable set and $U^\nu,U^\infty : E\larrow \mK$ are measurable maps, $\nu\in \N$. Then, there exist subsequences $(\nu_k)_1^\infty$, $(\nu_l)_1^\infty$:
\[
\begin{split}
& (1)\quad  \de_{U^\nu} \overset{_*}{\lharpoonup}\de_{U^\infty} \text{ in }\mY(E,\mK)  \ \ \Longrightarrow \ \ U^{\nu_k} \larrow U^\infty \ \text{ a.e.\ on }E.\\
& (2) \quad U^\nu \larrow U^\infty \ \text{ a.e.\ on }E\ \ \ \ \Longrightarrow \ \ \de_{U^{\nu_l}} \overset{_*}{\lharpoonup}\de_{U^\infty} \text{ in }\mY(E,\mK).
\end{split}
\]
\end{lemma}

\noi \textbf{General frames, derivative expansions, difference quotients.} In what follows we will consider non-standard orthonormal frames of $\R^{Nn^p}_s$ and write derivatives $D^p u$ with respect to them. This generalisation is irrelevant to the mollification results we establish herein but it was absolutely essential for the existence-uniqueness results we established in \cite{K8}-\cite{K10}. In any case, these bases will not appear explicitly anywhere in the proofs and they will not imply any technical ramifications. 

Let $\{E^1,...,E^N\}$ be an orthonormal frame of $\R^N$ and suppose that for each $\al=1,...,N$ we have an orthonormal frame $\{E^{(\al)1},...,E^{(\al)n}\}$ of $\R^n$. Given such bases, we will equip the space $\R^{Nn^p}_s$ with the following induced orthonormal base: 
\beq \label{2.2}
\ \ \R^{Nn^p}_s  =\, \spn[ \big\{ E^{\al i_1... i_p}  \big\}], \ \ \ \ \ E^{\al i_1... i_p}:= \,E^\al \ot \left( E^{(\al)i_1} \vee ... \vee E^{(\al)i_p} \right)
\eeq
where
\beq \label{2.1}
a \vee b \, :=\, \frac{1}{2}\Big(a\ot b \, +\, b\ot a \Big), \ \ \ a,b\in \R^n,
\eeq
is the symmetrised tensor product. Given such frames, let 
\[
D^p_{ E^{(\al)i_p} ... E^{(\al)i_1} }\, =\, D_{E^{(\al)i_p}}\cdots D_{E^{(\al)i_1}}
\]
denote the usual $p$th order directional derivative along the respective directions. Then, the $p$th order derivative $D^pu$ of a map $u:\Om \sub \R^n \larrow \R^N$ can be expressed as
\beq  \label{2.3}
\begin{split}
D^pu\, &=\, \sum_{\al, i_1,...,i_p} \Big( E^{\al i_1...i_p} : D^p u \Big) \, E^{\al i_1...i_p} =\\
& =\, \sum_{\al, i_1,...,i_p}  \Big( D^p_{E^{(\al)i_1} ...E^{(\al)i_p}}(E^\al \cdot u)\Big) \, E^{\al i_1...i_p} .
\end{split}
\eeq
We will use the following compact notation for the (formal) Taylor expansion around a point $x\in \Om$:
\[
u(z)\ =\ \sum_{p=1}^\infty \, \frac{1}{p!} \, D^{p}u(x) : (z-x)^{\ot p} .
\]
The notation ``${\ot p}$'' stands for the $p$th tensor power and ``:" is the obvious contraction of indices which in index form reads
\[
u_\al (z)\ =\ \sum_{p=1}^\infty \, \frac{1}{p!} \sum_{i_1,...i_p=1}^n D_{i_1...i_p}^{p}u_\al (x) \, (z-x)_{i_1} ...  (z-x)_{i_p} .
\]
Expansions analogous to \eqref{2.3} will also be applied to difference quotients which play a crucial role in our approach. Given $a\in \R^n$ with $|a|=1$ and $h\in \R\set \{0\}$, the $1$st order difference quotient of $u$ along the direction $a$ at $x$ will be denoted by
\beq  \label{2.4}
D^{1,h}_a u(x)\, :=\, \frac{u(x+ha)-u(x)}{h}.
\eeq 
By iteration, if $h_1,...,h_p\neq 0$ the $p$th order difference quotient along $a_1,...,a_p$ is
\beq  \label{2.5}
D^{p,h_p...h_1}_{a_p...a_1} u\, :=\, D^{1,h_p}_{a_p}\Big( \cdots \big( D^{1,h_1}_{a_1}u\big) \Big).
\eeq

We now introduce difference quotients taken with respect to frames as in \eqref{2.2}.

\begin{definition}[Difference quotients] 
\label{definition6} Let $\{E^1,...,E^N\}$ be an orthonormal frame of $\R^N$ and let also $\{E^{(\al)1},...,E^{(\al)n}\}$ be for each $\al=1,...,N$ an orthonormal frame of $\R^n$, while for any $p\in \N$ the tensor space $\R^{Nn^p}_s$ is equipped with the frame \eqref{2.2}. Given any \textbf{vector-indexed} infinitesimal sequence 
\[
\begin{split}
 (h_{\overline{m}} )_{\overline{m} \in \N^p} \sub \big(\R\set\{0\}\big)^{p},\ \ \ \overline{m}=( {m^1},...,{m^p} ), \ \ h_{m^q} \ri 0\ \text{ as } m^q\ri\infty,
\end{split}
\]
we define the \textit{$p$th order difference quotients} of the  measurable mapping $u:\Om\sub \R^n \larrow \R^N$ (with respect to the fixed reference frames) arising from $(h_{\overline{m}} )_{\overline{m} \in \N^p}$ as the family of maps 
\[
 D^{p,h_{\overline{m}}}u \ : \ \ \Om \sub \R^n \larrow \R^{Nn^p}_s,  \ \ \  \ \overline{m}\in \N^p,
 \]
each of which is given by
\[
\begin{split}
 D^{p,h_{\overline{m}}}u  \,  := \sum_{\al, i_1,...,i_p}  \Big[ D^{p,h_{m^p}...h_{m^1} }_{E^{(\al)i_p} ...E^{(\al)i_1}}(E^\al \cdot u)\Big]  \, E^{\al i_1...i_p} .
\end{split}
\]
The notation in the bracket above is as in \eqref{2.4}, \eqref{2.5}. Further, given any \textbf{matrix-indexed} infinitesimal sequence
\[
 (h_{\underline{m}} )_{\underline{m} \in \N^{p\by p}} \sub \big(\R\set\{0\}\big)^{p\by p}, \ \ \ \underline{m}  = 
\left[
\begin{array}{ccc}
m^1_1 & 0          & 0 \ \, ...\, \ 0 \\
m_2^1 & m_2^2 & 0 \ \, ... \, \ 0 \\
\vdots  &  & \ddots  \ \, \vdots \\
m_p^1 & m_p^2 & \,  ...\, \ \  m_p^p
\end{array}\right] \!, \ \
h_{m_p^q}\ri0\text{ as }m_p^q\ri\infty,
\]
we will denote its nonzero row elements by
\[
\underline{m}_q\, :=\, (m_q^1,...,m_q^q) \ \in \N^q,\quad q\,=\, 1,...,p,
\]
and we define the \textbf{$p$th order Jet $D^{[p],h_{\underline{m}}}u$ of difference quotients of $u$} (with respect to the reference frames) arising from $ (h_{\underline{m}} )_{\underline{m} \in \N^{p\by p}} $ as the family of maps
\[
 D^{[p],h_{\underline{m}}}u \ : \ \ \Om \sub \R^n \larrow \R^{Nn}\!\by ...\by \R^{Nn^p}_s, \ \ \  \ \underline{m}\in \N^{p\by p},
 \]
each of which is given by
\[
D^{[p],h_{\underline{m}}}u  \, :=\, \left(D^{1,h_{\underline{m}_1}}u,... \, ,D^{p,h_{\underline{m}_p}}u\right) .
\]
\end{definition}

\begin{definition}[Multi-indexed convergence] \label{Multi-indexed convergence} If $\underline{m}\in \N^{p\by p}$ is a lower trigonal matrix of indices as above, the expression ``$\underline{m}\ri \infty$" will by definition mean successive convergence with respect to each index separately in the following order: 
\[
\lim_{\underline{m}\ri \infty}\  :=\, \lim_{m^p_p \ri \infty}\, ...\,  \lim_{m^2_2 \ri \infty} \,\lim_{m_2^1 \ri \infty} \, \lim_{m^1_1 \ri \infty}.
\]
\end{definition}

\subsection{Main definitions and some analytic properties}    

\begin{definition}[Diffuse Jets]  \label{Diffuse Derivatives}
Suppose we have fixed some reference frames as in Definition \ref{definition6}. For any measurable map $u : \Om\sub \R^n\larrow\R^N$, we define the \textbf{diffuse $p$th order Jets $\mD^{[p]} u$ of $u$} as the following subsequential weak* limits:
\[
\de_{D^{[p],h_{\underline{m}}}u} \, \overset{^*}{\smash{\lharpoonup}}\ \mD^{[p]} u \ \ \   \text{ in }\mY\Big(\Om, \smash{\overline{\R}}^{Nn}\! \by ...\by \smash{\overline{\R}}^{Nn^p}_s\Big),
\]
which arise as $\underline{m}\ri \infty$ along multi-indexed infinitesimal subsequences.
\end{definition} 

As a consequence of the separate convergence, the $p$th order Jet is always a (fibre) product Young measure: $\mD^{[p]}u= \mD u \by \cdots \by \mD^p u$. 

Next is the central notion of generalised solution. We will use the notation $\underline{\X} \, \equiv \, (\X_1,...,\X_p)$ for points in $\smash{\overline{\R}}^{Nn} \by\cdots \by \smash{\overline{\R}}^{Nn^p}_s$ and also the symbol ``$\supp_*$" to denote the \emph{reduced support} of a probability measure $\vartheta \in \mP\big( \smash{\overline{\R}}^{Nn}\! \by ...\by \smash{\overline{\R}}^{Nn^p}_s\big)$ off ``infinity":
\[
\supp_* (\vartheta )\, :=\, \supp  ( \vartheta ) \cap \Big(\R^{Nn}\! \by \cdots \by \R^{Nn^p}_s \Big).
\]

\begin{definition}[$\mD$-solutions for $p$th order systems] \label{definition13} Let $\Om\sub \R^n$ be an open set and
\[
F\ : \ \Om \by \Big(\R^N\by \R^{Nn} \by \cdots \by \R^{Nn^p}_s \Big) \larrow \R^M
\]
a Carath\'eodory mapping. Assume also that we have fixed some reference frames as in Definition \ref{definition6} and consider the $p$th order PDE system
\beq \label{2.11a}
F\left(x,u(x),D^{[p]}u(x)\right)\, =\,0, \ \ \ x\in \Om.
\eeq
We say that the measurable map $u : \Om\sub \R^n \larrow \R^N$ is a \textbf{$\mD$-solution of \eqref{2.11a}} when for any diffuse $p$th order  Jet $\mD^{[p]}u$ of $u$  arising from any  infinitesimal multi-indexed sequence $(h_{\underline{m}})_{\underline{m} \in \N^{p\by p}}$ (Definition \ref{Diffuse Derivatives}) we have
\[
\sup_{\underline{\X}  \in \, \supp_* (\mD^{[p]}u(x)) }\,
\left| F\big( x,u(x),\underline{\X}\big)\right| \, =\, 0, \quad \text{ a.e. }x\in\Om.
\]
\end{definition}
We now consider the consistency of the $\mD$-notions with the strong/classical notions of solution. For more details we refer to \cite{K8} and also to \cite{K9,K10}. In general, diffuse derivatives may be nonunique for nonsmooth maps. However, as the next simple consequence of Lemma \ref{lemma2} shows, they are compatible with weak derivatives and a fortiori with classical derivatives:

\begin{lemma}[Compatibility of weak and diffuse derivatives] \label{lemma8} If $u \in W^{p,1}_{\text{loc}}(\Om,\R^N)$, then the $p$th order diffuse Jet $\mD^{[p]}u$  is unique and  for any $k\in \N$ we have
\[
\mD^{[p+k]}u\, =\, \de_{(Du,...,D^{p}u)}\by \mD^{p+1} \by ... \by \mD^{p+k} u, \ \ \text{ a.e.\ on } \Om.
\]
\end{lemma}

The next result asserts the plausible fact that $\mD$-solutions are compatible with strong solutions. Its proof is an immediate consequence of  Lemma \ref{lemma8}.

\begin{proposition}[Compatibility of strong and $\mD$-solutions] \label{proposition9} Let $F$ be a Carath\'eodory map as in \eqref{1.1} and $u\in W^{p,1}_{\text{loc}}(\Om,\R^N)$. Consider the $p$th order PDE system
\[
F\left(x,u(x) ,D^{[p]}u(x)\right)\, =\,0, \ \ x\in \Om.
\]
Then, $u $ is a $\mD$-solution on $\Om$ if and only if $u$ is a strong a.e.\ solution on $\Om$.
\end{proposition}

Lemma \ref{lemma8} and Proposition \ref{proposition9} remain true if $u$ is merely $p$-times differentiable in measure, a notion weaker than approximate differentiability (see \cite{K8,AM}). For more details on the material of this section (e.g.\ analytic properties, equivalent formulations of Definition \ref{definition13}, etc) we refer to \cite{K8}-\cite{K10}.

\section{Mollification of $\mD$-solutions to fully nonlinear systems} \label{section3}

We begin with the next result which is the main technical core of our constructions. Our method of proof is inspired by the paper of Alberti \cite{A}.

\begin{lemma}[Construction of the approximations] \label{lemma10} Let $\Om\sub \R^n \larrow \R^N$ be a measurable map and $p\in \N$. Then, for any $\e>0$ and any  multi-index $\underline{m}\in \N^{p\by p}$ as in Definition \ref{definition6} , there exist a measurable set $E_{\e,\underline{m}} \sub \Om$ and a smooth map $u^{\e,\underline{m}} \in C^\infty_0(\Om,\R^N)$ such that
\beq \label{3.1}
\left\{\ \ \ 
\begin{split}
& \big| E_{\e,\underline{m}} \big| \, \leq \, \e, \ms\\
& \big\| u \,-\, u^{\e,\underline{m}} \big\|_{L^\infty(\Om\set E_{\e,\underline{m}} )} \, \leq \, \e,\ms\\
& \left\|  D^{[p],h_{\underline{m}}}u \,-\, D^{[p]} u^{\e,\underline{m}} \right\|_{L^\infty(\Om\set E_{\e,\underline{m}} )} \,  \leq \, \e,
\end{split}
\right.
\eeq 
where $ D^{[p],h_{\underline{m}}}u$ is the $p$th order Jet of difference quotients of $u$. If moreover $|\Om|<\infty$, then $u^{\e,\underline{m}} \in C^\infty_c(\Om,\R^N)$.  Finally, if $u\in L^r(\Om,\R^N)$ for some $r\in[1,\infty)$, then we also have
\beq \label{3.2}
\big\| u - u^{\e,\underline{m}} \big\|_{L^r(\Om )} \, \leq \, \e.
\eeq
\end{lemma}

The reader can easily be convinced that even if $u\in L_{\text{loc}}^1(\Om,\R^N)$, the standard mollifier $u*\eta^\e$ of $u$ does not satisfy these approximation properties (the best we can get is approximation in dual spaces, not almost uniform on $\Om$).

\BPL \ref{lemma10}. Let $u:\Om\sub \R^n \larrow \R^N$ be a given measurable map (extended on $\R^n\set \Om$ by zero). Let $D^{[p],h_{\um}}u$ be the Jet of $p$th order difference quotients of $u$ where the multi-index $\underline{m}\in \N^{p\by p}$ is fixed. We also fix $\e>0$.

\ms

\noi \textbf{Step 1.} We may assume that $\Om$ has finite measure. This hypothesis does not harm generality for the following reason: assuming we have established \eqref{3.1}, \eqref{3.2}  on subdomains of $\Om$ which have finite measure, we can fill $\Om$ a.e.\ by disjoint open cubes $(\Om_i)_1^\infty$ such that $\big| \Om \set \big(\cup_1^\infty \Om_i\big) \big|= 0
$ and on each $\Om_i$ \eqref{3.1} holds with $2^{-i}\e$ instead of $\e$ for respective sequences of functions $(u^{\e,\um,i})_{i=1}^\infty$ and sets $(E_{\e,\um,i})_{i=1}^\infty$. Then, we define $u^{\e,\um} \in C^\infty_0(\Om,\R^N)$ and $E_{\e,\um}\sub \Om$ with $|E_{\e,\um}|\leq \e$ by taking 
\[
u^{\e,\um}|_{\Om_i}\, :=\, u^{\e,\um,i}, \ \ \ \  E_{\e,\um}\, :=\, \cup_{i=1}^\infty E_{\e,\um,i}. 
\]
The conclusion of Lemma \ref{lemma10} then follows. 

\ms

\noi \textbf{Step 2.} We now show there exists a measurable set and smooth maps 
\[
F_{\e,\um} \,\sub \, \Om,\ \ \ \ \big\{U^{q,\e,\um}\big\}_{q=0}^p, \ \  U^{q,\e,\um}\, \in C^\infty_c\big(\Om,\R^{Nn^q}_s\big),
\]
such that
\beq \label{3.4}
\left\{\ \ \ 
\begin{split}
& \big| F_{\e,\underline{m}} \big| \, \leq \, \e, \ms\\
& \big\| u \,-\, U^{0,\e,\underline{m}} \big\|_{L^\infty(\Om\set F_{\e,\underline{m}} )} \, \leq \, \e,\ms\\
& \left\|  D^{q,h_{\underline{m}}}u \,-\, U^{q,\e,\underline{m}} \right\|_{L^\infty(\Om\set F_{\e,\underline{m}} )} \,  \leq \, \e, \ \ q=1,...,p.
\end{split}
\right.
\eeq 
Indeed, let us define
\beq  \label{3.5}
V\,:=\, \big(u,D^{[p],h_\um }u \big)\ : \ \ \Om \sub \R^n \larrow \R^N\by\R^{Nn}\by\cdots\by\R^{Nn^p}_s
\eeq
and for any $R>0$ we consider the truncation
\beq \label{3.3}
T^R(\xi)\, := \left\{
\begin{split}
\xi,\ \ \ \ \ \text{ if } |\xi|<R, \\
R\frac{\xi}{|\xi|}, \ \text{ if }|\xi|\geq R.\\
\end{split}
\right.
\eeq
Then, we have that $|T^R(V)|\in (L^1\cap L^\infty)(\Om)$. Since $T^R(V) \larrow V$ a.e.\ on $\Om$ as $R\ri \infty$, we have that $T^R(V) \larrow V$ in measure on $\Om$ as $R\ri \infty$. By the identity
\[
\big|T^R(V)-V \big|\, =\, \big(|V|-R\big)\chi_{\{|V|\geq R\}},
\]
for any $R> \e$ we have that
\[
\begin{split}
\left|\big\{ T^{R}(V)\neq  V \big\}\right| \,&= \, \big|\big\{|V|\geq R \big\} \big|\\
&\leq \, \Big|\big\{ |V|\geq R/2 \big\} \cap \big\{ |V|> (R+\e)/2 \big\}\Big|\\
& =\, \Big|\Big\{ \big(|V|-R/2\big)\chi_{\{|V|\geq R/2\} } >\e/2 \Big\}\Big|\\
& =\, \Big|\Big\{\big|T^{R/2}(V)-V \big|>\e/2 \Big\}\Big|\\
& \!\! \larrow \, 0,
\end{split}
\]
as $R\ri\infty$. Hence, for $R(\e)>0$ large enough, there is a measurable set $A_\e \sub \Om$ such that
\beq \label{3.6}
 |A_\e|\leq \e/2, \ \ \ \ T^{R(\e)}(V)\,=\, V \text{ on }\Om\set A_\e.
\eeq
Further, we can find a sequence of smooth compactly supported maps $(V^{\e,k})_{k=1}^\infty$ such that for any $s\in [1,\infty)$,
\[
\text{$\big|V^{\e,k}\,  -\, T^{R(\e)}(V)\big|\larrow 0$, \ \  in $L^s(\Om)$ and a.e.\ on $\Om$ as $k\ri \infty$.} 
\]
By Egoroff's theorem, the convergence is almost uniform on $\Om$ as well. Hence, we can find a large enough $k(\e)\in \N$ and a measurable set $B_\e \sub \Om$ such that
\beq \label{3.7}
|B_\e| \leq \e/2,\ \ \ \ \big\| V^{\e,k(\e)} -\, T^{R(\e)}(V) \big\|_{L^\infty(\Om\set B_\e)} \leq\, \e.
\eeq
By \eqref{3.6} and \eqref{3.7}, we conclude that by considering the measurable set $F_{\e,\um}\sub \Om$ given by
\[
F_{\e,\um} \, :=\, A_\e \cup B_\e
\]
we have
\[
\big|F_{\e,\um}\big|\leq\e,\ \ \ \ \big\| V^{\e,k(\e)} -\, V \big\|_{L^\infty(\Om\set F_{\e,\um})} \leq\, \e.
\]
Hence, \eqref{3.4} ensues because in view of \eqref{3.5} we may define
\beq \label{3.8}
\left\{ \ \
\begin{split}
U^{0,\e,\um}\, &:=\, \text{Proj}_{\R^N} \big(V^{\e,k(\e)}\big),\\
 U^{q,\e,\um}\, &:=\, \text{Proj}_{\R^{Nn^q}_s}\big(V^{\e,k(\e)}\big),\ \ \ q=1,...p.\\
\end{split}
\right.
\eeq

\noi \textbf{Step 3.} We now establish \eqref{3.1} by using the previous step. Let $\om_{\e,\um}\in C^0[0,\infty)$ be the joint modulus of continuity of the mappings $\{U^{q,\e,\um}\}_{q=0}^p$ of \eqref{3.4}, i.e.\ $\om_{\e,\um}(0)=0$, $\om_{\e,\um}$ is increasing and 
\beq \label{3.9}
\sum_{q=0,1,...,p}\Big|U^{q,\e,\um}(z) \, -\, U^{q,\e,\um}(y)\Big| \, \leq\, \om_{\e,\um}\big(|z-y|\big),\ \ \ z,y\in \Om.
\eeq
We fix $\de \in (0,1)$ that will be specified later and consider the grid $(\de \N)^n\sub \R^n$. Let also
\beq \label{3.10}
\big\{ Q_{\de,i}\big\}_{i=1}^\infty \, :=\, \Big\{ \text{enumeration of open cubes whose vertices are on }(\de \N)^n \Big\} .
\eeq
Fix also $\al \in (0,1)$ and consider 
\beq \label{3.11}
\big\{ Q_{\al\de,i}\big\}_{i=1}^\infty \, :=\, \Big\{ \forall\, i,\text{ $Q_{\al\de,i}$ is concentric cube to $Q_{\de,i}$ with side length $\al\de$} \Big\}.  
\eeq
We further define
\beq \label{3.12}
\Om_\de\, := \bigcup_{i\in \N: \, Q_{\de,i}\sub \Om} Q_{\de,i}, \ \ \ \ \ \Om_{\al\de}\, :=\bigcup_{i\in \N: \, Q_{\de,i}\sub \Om} Q_{\al\de,i}.
\eeq
We now show that we may select $\al,\de>0$ such that
\begin{align}
\label{3.13} \max_{k=0,1,...,p} \Bigg\{ \sum_{q=k+1} & \frac{1}{(q-k)!}\, \big\| U^{q,\e,\um}\big\|_{L^\infty(\Om)} \,\de^{q-k} \Bigg\}\, \leq\, \e\\
\label{3.14} &\ \ \ \om_{\e,\um}(\de)\,\leq\, \e,\\
\label{3.15} & \ \big| \Om\set \Om_{\al\de} \big|\, \leq\, \e.
\end{align}
Indeed, since $\om_{\e,\um}$ is increasing and in view of \eqref{3.10}-\eqref{3.12}, by choosing $\de$ small enough we immediately obtain \eqref{3.13}-\eqref{3.14} and also 
\beq
\label{3.16}  \big| \Om\set \Om_{\de} \big|\, \leq\, \e.
\eeq
In view of the identity
\beq \label{3.16a} 
|\Om_{\al\de}|\, =\, \al^n\, |\Om_\de|,
\eeq
by choosing 
\beq \label{3.16b} 
\al\, :=\, \left( \max\left\{ 1- \frac{\e}{2|\Om|},\, \frac{1}{2} \right\}\right)^{1/n},
\eeq
we obtain
\[
\begin{split}
\big|\Om \set \Om_{\al\de}\big|\, &=\ \big|\Om \set \Om_{\de}\big|\, +\, \big|\Om_\de \set \Om_{\al\de}\big|\\
&\leq\ \e/2\, +\, |\Om_\de|\,-\,  \al^n |\Om_\de| \\
&\leq\ \e/2\, +\, (1\,-\,  \al^n) |\Om| \\
& \leq\ \e
\end{split}
\]
and \eqref{3.15} has been established as well. For each $i\in \N$ such that $Q_{\de,i}\sub \Om$, we consider a cut off function $\zeta^{\de,i} \in C^\infty_c(Q_{\de,i})$ such that
\beq \label{3.17}
\chi_{Q_{\al\de,i}}\, \leq\,\zeta^{\de,i}\, \leq\, \chi_{Q_{\de,i}}
\eeq
and let 
\[
\big\{x^{\de,i}\big\}_1^\infty :=\, \Big\{ \text{the centres of the cubes }\big\{ Q_{\de,i}\big\}_{i=1}^\infty\Big\}\sub \R^n .
\]
We now define a mapping 
\[
u^{\e,\um} \, \in \, C^\infty_c(\Om,\R^N)
\]
which is given by
\beq \label{3.18}
u^{\e,\um}(x)\, :=\, \sum_{i\in \N:\, Q_{\de,i}\sub \Om}\zeta^{\de,i}(x)\Bigg( \, \sum_{q=0}^p\, \frac{1}{q!}\, U^{q,\e,\um}\big(x^{\de,i}\big) : \big( x-x^{\de,i} \big)^{\ot q} \Bigg).
\eeq
Then, by \eqref{3.17}, on each $Q_{\al\de,i} \sub \Om$, the map $u^{\e,\um}$ equals the restriction of a $p$th order polynomial and hence for any $k\in\{1,...,p\}$, we have
\[
D^ku^{\e,\um}(x)\, =\,  \sum_{q=k}^p\, \frac{1}{(q-k)!}\, U^{q,\e,\um}\big(x^{\de,i}\big) : \big( x-x^{\de,i} \big)^{\ot (q-k)}  ,
\]
for $x\in Q_{\al\de,i}$. By recalling the properties of the measurable set $F_{\e,\um}\sub \Om$ of \eqref{3.4}, for any $k\in\{0,1,...,p\}$ and for a.e.\ $x\in Q_{\al\de,i}\set F_{\e,\um}$, we have
\[
\begin{split}
\Big| D^{k,h_\um}u\, -\, D^ku^{\e,\um}  \Big|(x)\, &\leq\, \Big| D^{k,h_\um}u   \, -\, U^{k,\e,\um} \Big|(x)\, +\, \Big| U^{k,\e,\um} \, -\, D^ku^{\e,\um} \Big|(x)\\
& \leq\ \e \, +\, 
\Bigg|\, U^{k,\e,\um}(x)\ -\\
&\ \ \ \ \ \ \  -\, \sum_{q=k}^p\, \frac{1}{(q-k)!}\, U^{q,\e,\um}\big(x^{\de,i}\big) : \big( x-x^{\de,i} \big)^{\ot (q-k)} \Bigg|
\end{split}
\]
which by \eqref{3.9} gives that 
\[
\begin{split}
\Big| D^ku^{\e,\um}(x)  \, -\,D^{k,h_\um}u (x)\Big|\, &\leq\ \e \, +\,    \Big|  U^{k,\e,\um}(x) \, -\, U^{k,\e,\um}(x^{\de,i} ) \Big|  \\
&\ \ \ \ \ \ \ +\,     \sum_{q=k+1}^p\, \frac{1}{(q-k)!}\, \big\| U^{q,\e,\um}\big\|_{L^\infty(\Om)} \big| x-x^{\de,i} \big|^{q-k} \\
&\leq\ \e \, +\,  \om_{\e,\um}\big(\big| x-x^{\de,i} \big|\big)  \\
&\ \ \ \ \ \ \ +\,    \sum_{q=k+1}^p\, \frac{1}{(q-k)!}\, \big\| U^{q,\e,\um}\big\|_{L^\infty(\Om)} \, \de^{q-k}
\end{split}
\]
for a.e.\ $x\in Q_{\al\de,i}\set F_{\e,\um}$. Hence, 
\beq \label{3.19}
\begin{split}
\Big| D^ku^{\e,\um} \, -\, D^{k,h_\um}u \Big|\, \leq & \ \e \,+\,  \om_{\e,\um}(\de)  \\
& \ \ \ +\,       \sum_{q=k+1}^p\, \frac{1}{(q-k)!}\, \big\| U^{q,\e,\um}\big\|_{L^\infty(\Om)} \, \de^{q-k},
\end{split}
\eeq
a.e.\ on $\in Q_{\al\de,i}\set F_{\e,\um}$. By \eqref{3.19} and \eqref{3.10}-\eqref{3.14} we deduce that
\beq \label{3.20}
\max_{k=0,1,...,p} \Big\| D^ku^{\e,\um} \, -\, D^{k,h_\um}u \Big\|_{L^\infty(\Om_{\al\de}\set F_{\e,\um})}\, \leq\ 3\e. 
\eeq
Finally, we set
\beq \label{3.21}
E_{\e,\um}\, :=\, F_{\e,\um} \cup (\Om\set \Om_{\al\de}\big).
\eeq
Then, by \eqref{3.21}, \eqref{3.15} and \eqref{3.4} we have that
\[
\big| E_{\e,\um} \big| \, \leq \,2\e
\]
and  \eqref{3.20}, \eqref{3.21} imply
\beq \label{3.22}
\left\{
\begin{split}
 \Big\| D^{[p]}u^{\e,\um} \, -\, D^{{[p]},h_\um}u \Big\|_{L^\infty(\Om \set E_{\e,\um})}\, &\leq\ 3p\,\e,\\
 \big\| u^{\e,\um} \, -\, u \big\|_{L^\infty(\Om\set E_{\e,\um})}\, &\leq\ 3\e.
 \end{split}
 \right.
\eeq
Hence, by replacing $\e$ by $\e/3p$, we see that \eqref{3.1} has been established.

\ms

\noi \textbf{Step 4.} We now establish \eqref{3.2} under the additional hypothesis that $u\in L^r(\Om,\R^N)$. We begin by noting that on top of \eqref{3.4} we can also arrange to have
\beq  \label{3.22}
\big\| u\, -\, U^{0,\e,\um}\big\|_{L^r(\Om)}\, \leq\, \e.
\eeq
This follows by the next simple modification of Step 2: we replace $T^R(V)$ by $\big(u,T^R\big(D^{[p],h_\um}u \big)\big)$, choose $s:=r$ and use that $u$ can be approximated in the $L^r$ norm by smooth compactly supported mappings. Then we obtain \eqref{3.22} and the first and last inequalities of \eqref{3.4}. The middle inequality of \eqref{3.4} follows by \eqref{3.22} and (by perhaps modifying $F_{\e,\um}$ and the choice of $\e$ accordingly):
\[
\begin{split}
\Big| \Big\{ \big|u-U^{0,\e,\um}\big|>\e^{1/2} \Big\}\Big| \, &\leq\,  \e^{-r/2}\int_\Om \big|u\,-\,U^{0,\e,\um}\big|^r \\
 & \leq\, \e^{r/2}.
 \end{split}
\]
Next, by \eqref{3.22} and by \eqref{3.18} we have that
\[
\begin{split}
\big\|u-u^{\e,\um} \big\|^r_{L^r(\Om)} \, &\leq\, 2^r\e\, +\, 2^r\int_{\Om}\Big|u^{\e,\um} \, -\, U^{0,\e,\um}\Big|^r\\
&\leq\, 2^r\e\, +\, 2^r\big|\Om\set\Om_\de \big| \,   \big\|U^{0,\e,\um}\big\|_{L^\infty(\Om)}^r\\
&\ \ \ \ \ \ \ \ \  +\, 
4^r\int_{\Om_\de \set\Om_{\al\de}}\Big\{ \big|u^{\e,\um}|^r+\big|U^{0,\e,\um}\big|^r\Big\}\\
&\ \ \ \ \ \ \ \ \  + \, 2^r\int_{\Om_{\al\de}}\Big|u^{\e,\um}\, -\, U^{0,\e,\um}\Big|^r
\end{split}
\]
which gives
\begin{align}  \label{3.24}
\big\|u-u^{\e,\um} \big\|^r_{L^r(\Om)} \, &\leq\, 4^r\e\, +\, 4^r\big|\Om\set\Om_\de \big|\,  \big\|U^{0,\e,\um}\big\|_{L^\infty(\Om)}^r \nonumber\\
&\ \ \ \ \ \ \ \ \ \, +\, 4^r \sum_{i\in\N:\, Q_{\de,i}\sub \Om} \,\Bigg\{ \int_{Q_{\de,i} \set Q_{\al\de,i}}\Big\{\big|u^{\e,\um}|^r+\big|U^{0,\e,\um}\big|^r\Big\}\\
&\hspace{90pt}\ \ \ \ \, + \int_{Q_{\al\de,i}}\Big|u^{\e,\um}\, -\, U^{0,\e,\um}\Big|^r\Bigg\}. \nonumber
\end{align}
Moreover, by \eqref{3.18} (and \eqref{3.17}) we have the estimates
\beq  \label{3.25}
\left\{ \ \,
\begin{split}
\big| u^{\e,\um}\big|\, &\leq\, \sum_{q=0}^p \,\frac{1}{q!}\, \big\|U^{q,\e,\um}\big\|_{L^\infty(\Om)} \,\de^q ,\ \ \ \text{ on }Q_{\de,i}\set Q_{\al\de,i},\ i\in\N,\\
\big| u^{\e,\um} \,-\,U^{0,\e,\um}\big|\, &\leq\, \sum_{q=1}^p \,\frac{1}{q!}\, \big\|U^{q,\e,\um}\big\|_{L^\infty(\Om)} \,\de^q ,\ \ \ \text{ on }Q_{\al\de,i},\ i\in\N.
\end{split} 
\right.
\eeq
By inserting \eqref{3.25} into \eqref{3.24} we obtain the estimate
\[
\begin{split}
\big\|u\,-\,u^{\e,\um} \big\|^r_{L^r(\Om)} \, &\leq\, 4^r\e\, +\, 4^r \big|\Om\set\Om_\de \big|\,  \big\|U^{0,\e,\um}\big\|_{L^\infty(\Om)}^r\\
&\ \ \ \ \ \ \ \ \ \, +\, 6^{r}\big|\Om_\de \set\Om_{\al\de} \big|  \left(\, \sum_{q=0}^p \,\frac{1}{q!}\, \big\|U^{q,\e,\um}\big\|_{L^\infty(\Om)} \,\de^q\right)^r
\\
&\ \ \ \ \ \ \ \ \ \, +\, 4^r|\Om|  \left( \, \sum_{q=1}^p \,\frac{1}{q!}\, \big\|U^{q,\e,\um}\big\|_{L^\infty(\Om)} \,\de^q \right)^r
\end{split}
\]
which in turn gives
\beq \label{3.26}
\begin{split}
\big\|u\,-\,u^{\e,\um} \big\|_{L^r(\Om)} \, &\leq\, 4\e\, +\, 4 \big|\Om\set\Om_\de \big|^{1/r}\,  \big\|U^{0,\e,\um}\big\|_{L^\infty(\Om)}\\
&\ \ \ \ \ \ \ \ \, +\, 6(1-\al^n)^{1/r} \left( |\Om|^{1/r}  \sum_{q=0}^p \,\frac{1}{q!}\, \big\|U^{q,\e,\um}\big\|_{L^\infty(\Om)}   \right)
\\
&\ \ \ \ \ \ \ \ \,  +\, 4|\Om|  \, \sum_{q=1}^p \,\frac{1}{q!}\, \big\|U^{q,\e,\um}\big\|_{L^\infty(\Om)} \,\de^q  .
\end{split}
\eeq
In view of \eqref{3.26} and of \eqref{3.13}-\eqref{3.16b}, by decreasing $\de$ further and by increasing $\al$ even further, we can achieve
\[
\big\|u\,-\,u^{\e,\um} \big\|_{L^r(\Om)}\, \leq\,7\e.
\]
Hence, the desired conclusion follows by replacing $\e$ by $\e/7$.  The lemma has been established.       \qed

\ms

By utilising Lemma \ref{lemma10}, we may now state and prove the main result of this paper.

\begin{theorem}[Mollification of $\mD$-solutions] \label{theorem11} Let $\Om\sub \R^n$ be an open set and
\[
F\ : \ \Om \by \Big(\R^N\by \R^{Nn} \by \cdots \by \R^{Nn^p}_s \Big) \larrow \R^M
\]
a Carath\'eodory map. Suppose that $u : \Om \sub \R^n \larrow \R^N$ is a $\mD$-solution to the $p$th order PDE system
\[
F\Big(x,u(x),D^{[p]}u(x) \Big) \, =\, 0, \ \ \ \ x\in \Om.
\]
Then, there exists a multi-indexed sequence of maps $(u^{\underline{m}})_{\underline{m} \in \N^{p\by p}} \sub C^\infty_0(\Om,\R^N)$ with the following properties: for any diffuse $p$th order jet $\mD^{[p]} u$ generated along subsequences of an infinitesimal multi-indexed sequence $(h_{\underline{m}})_{\underline{m} \in \N^{p\by p}}$, there exists a single-indexed subsequence
\beq  \label{3.27}
(u^\nu)_{\nu \in \N} \, \sub \, (u^{\underline{m}})_{\underline{m} \in \N^{p\by p}} 
\eeq
with the following properties:
\beq  \label{3.28}
\left\{\ \ \ 
\begin{split}
& u^\nu \larrow u \hspace{45pt} \text{ a.e. on }\Om,\\
& \de_{D^{[p]}u^\nu} \, \overset{^*}{\smash{\lharpoonup}}\ \mD^{[p]} u\ \ \   \text{ in }\, \mY\Big(\Om, \smash{\overline{\R}}^{Nn} \by ...\by \smash{\overline{\R}}^{Nn^p}_s\Big),
\end{split}
\right.
\eeq 
both as $\nu \ri \infty$. In addition, for each $\nu\in\N$, the map $u^\nu \in C^\infty_0(\Om,\R^N)$ is a smooth strong solution to the approximate $p$th order PDE system
\beq  \label{3.29}
F\Big(x,u^\nu(x),D^{[p]}u^\nu(x) \Big) \, =\, f^\nu(x), \ \ \ \ x\in \Om
\eeq
and $f^\nu :\Om\sub \R^n \larrow\R^M$ is a measurable map satisfying $f^\nu \larrow 0$ as $\nu\ri \infty$ in the following sense:
\beq \label{3.30}
 \text{For any }\Phi \in C^0_c\big( \R^{Nn} \! \by \cdots \by \R^{Nn^p}_s \big),\ \ \ \Phi\big( D^{[p]}u^\nu\big) f^\nu \larrow 0, \text{ a.e. on }\Om.
 \eeq
\end{theorem}

Note that the second statement of \eqref{3.28} is interesting because the (fibre) product Young measure $\mD^{[p]} u$ can be weakly* approximated by the product Young measures $\de_{D^{[p]}u^\nu} $ as $\nu \ri\infty$ (confer with Definition \ref{Diffuse Derivatives}). The following consequence of our constructions will also follow from Theorem \ref{theorem11}.

\begin{corollary}[Mollification of $\mD$-solutions contn'd] \label{corollary12} In the setting of Theorem \ref{theorem11}, the mode of convergence \eqref{3.30} is (up to the passage to a subsequence) equivalent to either of the modes of convergence as $\nu\ri\infty$:
\beq  \label{3.31}
\left\{
\begin{split}
& (a) \text{ For any } R>0, \ \ \ \ \chi_{\mB_R(0)}\big( D^{[p]}u^\nu\big) f^\nu \larrow 0, \text{ a.e. on }\Om. \\
& (b)\text{ For any } \e>0, \text{ and any } E\sub \Om \text{ with } |E|<\infty,\\
&\hspace{65pt} \Big| E\cap\big\{ |f^\nu|>\e\big\} \cap \big\{ \big|D^{[p]}u^\nu\big|<{1}/{\e}\big\} \Big|\larrow 0. \ \ \ \ \ \ \ \ \ \ 
\end{split}
\right.
\eeq
In addition, we have
\beq  \label{3.32}
f^\nu \larrow 0 \ \ \ \text{ a.e. on }\Om\set E,
\eeq
where
\beq \label{3.33}
E\, :=\, \overset{p}{\underset{q=1}{\bigcup}}\Big\{x\in \Om\, :\, \big[\mD^{q}u(x)\big]\big(\{\infty\}\big) > 0 \Big\}.
\eeq
\end{corollary}

 The next example shows that in general it is not possible to strengthen the modes of convergence in \eqref{3.30}-\eqref{3.33} to $L^p$ or even to a.e.\ on $\Om$ even for linear equations:

\begin{example}[Optimality of Theorem \ref{theorem11}] Consider the equation 
\[
u'\, =\, 0, \ \ \ u \ : (0,1)\sub \R\ri \R. 
\]
(I) Let $u_c \in C^0_0(0,1)$ be a Cantor-type function\footnote{The 
reader should not be alarmed by the fact that $\mD$-solutions satisfy such ``unpleasant" nonuniqueness properties. On the one hand, this behaviour is a common feature of all ``a.e.-type" notions of solution, e.g.\ $L^p$-viscosity solutions, piecewise solutions, etc. On the other hand, in the vectorial nonlinear case, most interesting systems present such issues even for smooth $C^\infty$ solutions (e.g.\ the minimal surface system in codimension greater than $1$, see \cite{LO}). \textbf{However}, by imposing extra constraints on top of the boundary conditions (which depend on the nature of the problem), \textbf{$\mD$-solutions are indeed unique} and pathological solutions as in this example are ruled out, see \cite{K8} and 
\cite{K10}.}
with $u_c\neq 0$ and $u_c'=0$ a.e.\ on $(0,1)$. Then, $u_c$ is a $\mD$-solution because it satisfies the equation a.e.\ on $(0,1)$. However, for any hypothetical approximate equation 
\[
{u_c^\e}'\, =\, f^\e, \ \ \text{ on }(0,1)
\]
such that $u^\e_c \larrow u_c$ a.e.\ on $\Om$ and $f^\e \larrow 0$ in $L^1(0,1)$, by Poincar\'e inequality we have
\[
\| {u_c^\e}\|_{L^1(0,1)}\, \leq\, C \| {u_c^\e}'\|_{L^1(0,1)}\, =\, C \| f^\e\|_{L^1(0,1)} \larrow 0, \ \ \text{ as }\e\ri0,
\]
while $u^\e_c \larrow u_c$ a.e.\ on $\Om$ and $u_c\neq 0$.
\ms

\noi (II) Let $u_s \in L^1_c(0,1)$ be a singular solution to the equation such that for any diffuse derivative $\mD u_s \in \mY\big((0,1),\overline{\R}\big)$, we have $\mD u_s \equiv \de_{\{\infty\}}$ on a set of positive measure. Such a solution can be constructed by taking a compact nowhere dense set $K\sub (0,1)$ of positive measure (e.g., we may take
\[
K\, =\, [1/3,\,2/3] \set \cup_{j=1}^\infty\big( q_j-3^{-2j}, q_j+3^{-2j} \big)
\]
where $(q_j)_1^\infty$ is an enumeration of the rationals in $[1/3,\,2/3]$) and setting $u_s:=\chi_K$. Then, $u_s$ satisfies 
\[
\text{$|D^{1,h}u_s(x)| \larrow 0$ for $x\in (0,1)\set K$, \ \ \ $|D^{1,h}u_s(x)| \larrow \infty$ for $x\in K$,} 
\]
as $h\ri 0$. By Lemma \ref{lemma2}, we have that all diffuse gradients coincide and are given by
\[
\mD u_s(x) \, =\,  \chi_{ (0,1)\set K}(x)\, \de_{\{0\}}\, +\,  \chi_{K}(x) \, \de_{\{\infty\}}, \ \ \ \text{ a.e. }x\in (0,1).
\]
This implies that $u_s$ is a $\mD$-solution: indeed, we have that either $\supp_*(\mD u_s(x))=\{0\}$ or $\supp_*(\mD u_s(x))=\emptyset$ and hence
\[
\text{ for a.e. }x\in (0,1)\text{ and all }X_x\in \, \supp_*(\mD u_s(x)) \ \ \Longrightarrow \ \ |X_x| \, =\, 0.
\]
However, it is impossible to obtain $f^\e\larrow0$ a.e.\ on $(0,1)$ for any approximation scheme as in the theorem such that 
\[
{u_s^\e}'\, =\, f^\e\ \ \text{ on }(0,1), \ \ \ u^\e_s \larrow u_c \text{ a.e. on }\Om.
\]
Indeed, by \eqref{3.28} we would have $\de_{{u_s^\e}'}  \overset{^*}{\smash{\lharpoonup}}\, \mD u_s$  in $\mY\big((0,1),\overline{\R}\big)$ as $\e\ri 0$ and hence we have $\de_{{u_s^\e}'}  \overset{^*}{\smash{\lharpoonup}}\,  \de_{\{\infty\}}$ in $\mY\big(K,\overline{\R}\big)$. Hence, by  Lemma \ref{lemma2} this would give $ |f^\e(x)| \, =  \,  |{u_s^\e}' (x)| \larrow \infty$ for a.e.\ $x\in K$. 
\end{example}

For the proof of Theorem \ref{theorem11} we need three lemmas which are given right next. The first two are variants of results established in \cite{K8}, while the third one is a consequence of standard results on Young measures. They are all given below in the generality of Young measures because they do not utilise the special structure of diffuse derivatives. For the sake of completeness, we provide the first two results in full by giving all the details of their proofs, while for the third we give a precise reference. 

\begin{lemma}[Convergence lemma V2, cf.\ \cite{K8}]\label{lemma16A} Suppose that $u^m \larrow u^\infty$ a.e.\ on $\Om$, as $m \ri \infty$ where $u^\infty, (u^m)_1^\infty$ are measurable maps $ \Om\sub \R^n \larrow \R^N$. Let $\mathbb{W}$ be a finite dimensional metric vector space, isometrically and densely contained into a compactification $\mK$ of itself. Suppose also we are given Carath\'eodory mappings 
\[
F^\infty,\, F^m\ : \ \Om \by \big( \R^N \by \mathbb{W} \big) \larrow \R^M, \ \ \ m\in \N,
\]
and we are also given Young measures $\vartheta^\infty,(\vartheta)_1^\infty$ in $\mY\big( \Om,\mK\big)$ such that the following modes of convergence hold true:
\[
\begin{split}
F^\mu (x,\cdot , \cdot)  \larrow \, & F^\infty   (x , \cdot  , \cdot) \ \text{ in }C^0(\R^N\!\by \mathbb{W}), \ \ \text{as }m\ri \infty,\ \text{ a.e. }x\in\Om,
\\
& \vartheta^m \, \overset{^*}{\smash{\lharpoonup}}\  \vartheta^\infty \ \text{ in }\ \mY\big( \Om, \mK\big), \ \ \ \text{as }m\ri \infty.
\end{split}
\]
Then, if for a given function $\Phi \in C^0_c(\mathbb{W})$ we have
\[
\int_{\mK}    \big| \Phi(\X) \, F^\infty\big(x,u^\infty(x),\X\big)\big| \, d[\vartheta^\infty(x)](\X)\, =\, 0,  \ \ \text{ a.e. }x\in \Om,
\]
it follows that
\[
\lim_{m\ri \infty}\, \int_{\mK}  \big| \Phi(\X) \,  F^m \big(x,u^m (x),\X\big)\big| \,d[\vartheta^m(x)](\X)\, =\, 0,  \ \ \text{ a.e. }x\in \Om.
\]
\end{lemma}

We will later apply this lemma to the case of $\mathbb{W} =\R^{Nn}\!\by\cdots\by\R^{Nn^p}_s$ for the compactification of $\mathbb{W}$ given by the torus $\mK=\smash{\overline{\R}}^{Nn} \!\by\cdots \by \smash{\overline{\R}}^{Nn^p}_s$. We remind that the metric on $\mathbb{W}$ is the product metric induced by the imbedding of $\R^{Nn^q}_s$ into $\smash{\overline{\R}}^{Nn^p}_s$, $q=1,...,p$.

\BPL \ref{lemma16A}. We first fix $\Phi \in C^0_c(\mathbb{W})$ and define
\[
\phi^m(x)\, :=\, \Big\|\Phi(\cdot)\,  F^m\big(x,u^m(x),\cdot\big)\,-\, \Phi(\cdot)\, F^\infty\big(x,u^\infty(x),\cdot\big) \Big\|_{C^0(\mathbb{W})}
\] 
and we claim that our convergence hypotheses imply $\phi^m\larrow 0$ a.e.\ on $\Om$. Indeed, let us fix $x\in \Om$ such that $u^m(x)\larrow u^\infty(x)$ (and the set of such points $x$ has full measure). Then, we can find compact sets $K'\Subset \R^N$ and $K''\Subset \mathbb{W}$ such that for large $m\in \N$ we have $u^m(x),u^\infty(x)\in K'$ and also $\supp(\Phi)\sub K''$. By the convergence assumption on the maps $F^m$, we have
\[
\big\| F^m(x,\cdot)-F^m(x,\cdot) \big\|_{C^0(K' \by K'')}\larrow 0, \ \ \text{ as } m\ri \infty.
\]
If $\om^\infty_x \in C^0[0,\infty)$ is the modulus of continuity of the map $K' \ni \xi\mapsto F^\infty(x,\xi,\X) \in \R^M$ (that is $\om^\infty_x \geq \om^\infty_x(0)=0$) which is uniform with respect to $\X\in K''$,  we have
\[
\begin{split}
|\phi^m(x)|\, & \leq \,   \|\Phi\|_{C^0(K'')}\, \Bigg\{  \Big\| F^\infty \big( x,u^m(x),\cdot\big)\,-\,  F^\infty\big(x,u^\infty(x),\cdot\big) \Big\|_{C^0(K'')} \\
&\hspace{70pt} +\,   \Big\| F^m \big( x,u^m(x),\cdot\big)\,-\,  F^\infty\big(x,u^m(x),\cdot\big) \Big\|_{C^0(K'')}
 \Bigg\}\\
&\leq \,  \|\Phi\|_{C^0(\mathbb{W})} \Big\{\,  \om^\infty_x\big(\big|u^m(x)-u^\infty(x)\big| \big)\  +\, \big\| F^m(x,\cdot)-F^\infty(x,\cdot) \big\|_{C^0(K'\by K'')}\Big\}\\
& \!\!\! \larrow 0, 
\end{split}
\]
as $m\ri \infty$. Since this happens for a set of points $x\in \Om$ of full measure, we deduce that $\phi^m(x)\larrow 0$ for a.e.\ $x\in\Om$.  We now fix $R>0$ and $\Phi \in C^0_c(\mathbb{W})$ as in the statement of the lemma and set 
\[
{\Om_R}\, := \, \mB_R(0) \cap \Big\{ x\in \Om \, : \, \big\|\Phi(\cdot)\, F^\infty\big(x,u^\infty(x),\cdot\big) \big\|_{C^0(\mathbb{W})} <\, R   \Big\}. 
\]
Since $|{\Om_R}|<\infty$, by Egoroff's theorem we can find a measurable sequence $\{E_i\}_1^\infty \sub {\Om_R}$ such that $|E_i|\ri 0$ as $i\ri \infty$ and for each $i\in\N$ fixed,
\[
\|\phi^m\|_{L^\infty({\Om_R}\set E_i)}\larrow 0, \ \ \text{ as } m\ri \infty. 
\]
Since $|\Om_R|<\infty$, this gives $\phi^m\larrow 0 $ in $L^1({\Om_R}\set E_i)$ as $m\ri\infty$. This convergence and the form of $\Om_R$ imply that the Carath\'eodory functions
\[
\begin{split}
\Psi^m(x,\X)\, &:=\,\left|\Phi(\X)\, F^m\big(x,u^m(x),\X\big)\right|,\\
\Psi^\infty(x,\X)\, &:=\,\left|\Phi(\X)\, F^\infty\big(x,u^\infty(x),\X\big)\right| ,
\end{split}
\] 
are elements of the space $L^1\big({\Om_R}\set E_i, C^0(\mK)\big)$ because 
\[
\|\phi^m\|_{L^1({\Om_R}\set E_i)}\, \geq\, \| \Psi^m-\Psi\|_{L^1({\Om_R}\set E_i, C^0(\mK))}
\]
and for $m$ large we have
\[
\| \Psi^m\|_{L^1({\Om_R}\set E_i)}\, \leq\, 1\,+\, \| \Psi^\infty\|_{L^1({\Om_R}\set E_i)}\, \leq\, 1+\, |\Om_R| \, R.
\]
Moreover, we have 
\[
\Psi^m \larrow \Psi^\infty \ \ \text{ in }L^1\big({\Om_R}\set E_i, C^0(\mK)\big), \ \text{ as }m\ri\infty
\]
and also by assumption $\vartheta^m \overset{^*}{\smash{\lharpoonup}}\,  \vartheta^\infty$ in $\mY\big({\Om_R}\set E_i,\mK \big)$ as $m\ri \infty$. Thus, the weak*-strong continuity of the duality pairing
\[
 L^\infty_{w^*}\Big({\Om_R}\set E_i,\mM ( \mK) \Big) \, \by\,  L^1\Big( {\Om_R}\set E_i, C^0 (  \mK)\Big) \larrow \R,
\]
implies
\[
\begin{split}
\lim_{m\ri \infty}  \, \int_{ {\Om_R}\set E_i } & \int_{\mK}    \big|\Phi(\X) \,  F^m \big(x,u^m (x),\X\big)\big|\,  d[\vartheta^m(x)](\X)\,dx\\
 &=\,     \int_{ {\Om_R}\set E_i } \int_{\mK}  \big| \Phi(\X) \,F^\infty \big(x,u^\infty (x),\X\big) \big| \, d[\vartheta^\infty (x)](\X)\,dx.
 \end{split}
\]
We recall now that by our hypothesis the right hand side of the above vanishes in order to obtain the desired conclusion after letting $i\ri \infty$ and then taking $R\ri \infty$.            \qed

\ms

The next lemma says that if the distance between two sequences of measurable maps asymptotically vanishes, then the maps represent the same Young measure in the compactification (see also \cite{K8,FL}).

\begin{lemma} \label{lemma34} Let $\mathbb{W}$ be a finite dimensional metric vector space, isometrically and densely contained into a compactification $\mK$ of itself. Let also $E\sub \R^n$ be a measurable set. If $U^m,V^m : E\sub \R^n \larrow \mW$ are measurable maps satisfying
\[
\de_{U^m} \overset{^*}{\smash{\lharpoonup}}\, \,\vartheta\  \text { in }\mY(E,\mK),\ \ \ \ |U^m-V^m|  \larrow\, 0 \, \text{ a.e. on }E,
\]
as $m\ri \infty$, then $\de_{V^m}  \overset{^*}{\smash{\lharpoonup}}\, \, \vartheta$ in $\mY(E,\mK)$, as $m\ri\infty$.
\end{lemma}

\BPL \ref{lemma34}. We begin by fixing $\e>0$, $\phi \in L^1(E)$ and $\Phi \in C^0(\mK)$. Since $\Phi$ is uniformly continuous on $\mK$, there exists a bounded increasing modulus of continuity $\om \in C^0[0,\infty)$ such that 
\[ 
|\Phi(X) -\Phi(Y)| \, \leq \, \om(|X-Y|), \ \ \ \text{ for }X,Y \in \mK
\]
and we also have $\|\om\|_{C^0(0,\infty)}<\infty$. Moreover, since $|U^m-V^m|\larrow 0$ a.e.\ on $E$, we have that $|U^m-V^m| \larrow 0$ $\mu$-a.e.\ on $E$ where $\mu$ is the absolutely continuous finite measure given by $\mu(A):=\|\phi\|_{L^1(A)}$, $A\sub E$. It follows that $|U^m-V^m|  \larrow 0$ in $\mu$-measure as well. As a consequence, we obtain
\[
\begin{split}
\left| \int_E \phi \Big[  \Phi(V^m) - \Phi(U^m)\Big] \right|\, & \leq\, \int_E |\phi|\, \om \big(|U^m-V^m|\big)\\
 \leq \, \|\om\|_{C^0(0,\infty)}  & \,\mu\big(\{|U^m-V^m|>\e\}\big) \ +\ \om(\e)\, \mu(E).
\end{split}
\]
By letting first $m\ri \infty$ and then $\e\ri 0$, the density  in $L^1\big(E,C^0(\mK)\big)$ of the linear span of the products of the form $\phi(x)\Phi(X)$ implies the desired conclusion.          \qed

\ms

The following result is the last ingredient needed for the proof of Theorem \ref{theorem11}.

\begin{lemma} \label{lemma35} Let $\mW'$ and $\mW''$ be finite dimensional metric vector spaces, isometrically and densely contained into certain compactification $\mK'$ and $\mK''$ of $\mW'$, $\mW''$ respectively. Let also $E\sub \R^n$ be a measurable set. If $U^m : E\sub \R^n \larrow \mW'$, $V^m : E\sub \R^n \larrow \mW''$ are sequences of measurable maps satisfying
\[
U^m\larrow\, U \, \text{ a.e. on }E,\ \ \ \ \de_{V^m} \overset{^*}{\smash{\lharpoonup}}\, \,\vartheta\  \text { in }\mY(E,\mK''),
\]
as $m\ri \infty$, then 
\[
\ \ \text{ $\de_{(U^m,V^m)} \, \overset{^*}{\smash{\lharpoonup}}\ \de_U \by \vartheta$\ \ \ in $\mY(E,\,\mK'\!\by\mK'')$, \ as $m\ri\infty$.}
\]
\end{lemma}

\BPL \ref{lemma35}. The proof of this result can be found (actually in a much more general topological setting) e.g.\ in \cite{FG}, Corollary 3.89 on p.\ 257. \qed

\ms 

Now we may prove our main result.

\BPT \ref{theorem11}. \textbf{Step 1}. We begin by noting a general convergence fact. Let $(\mathfrak{X},\rho)$ be a metric space, $f\in \mathfrak{X}$ and let also $\um \in \N^{p\by p}$ be a matrix of indices as in Definition \ref{definition6}. Suppose that $\{f^\um\}_{\um \in \N^{p\by p}} \sub \mathfrak{X}$ is a multi-indexed sequence such that the successive limit converges to $f$ in $\mathfrak{X}$ (Definition \ref{Multi-indexed convergence}):
\[
\lim_{\um \ri\infty} \, \rho\left(f^\um ,\, f\right)
\,=\, 
\lim_{m^p_p \ri\infty} \left(... \, ...\lim_{m^2_2 \ri\infty}\, \left(\lim_{m^1_2 \ri\infty} \left( \lim_{m^1_1 \ri\infty} \, \rho\left(f^\um ,\, f\right)\right)\right)\right)\,=\,0.
\]
Then, there exist subsequences $(m^b_{a,\nu})_{\nu=1}^\infty$, $a,b\in\{1,...,p\}$, such that, if $(\um_\nu)_{\nu=1}^\infty \sub \N^{p\by p}$ is the single-indexed sequence with components $m^b_{a,\nu}$, then $f^{\um_\nu}$ converges to $f$ in $\mathfrak{X}$:
\[
\lim_{\nu \ri\infty} \, \rho \left(f^{\um_\nu} ,\, f\right) \,=\,0.
\]
This is a simple consequence of the definitions of limits.

\ms

\noi  \textbf{Step 2}. Let $u:\Om\sub\R^n \larrow \R^N$ be a $\mD$-solution to the system 
\[
F\Big(x,u(x),D^{[p]}u(x) \Big) \, =\, 0, \ \ \ \ x\in \Om,
\]
and let $\mD^{[p]} u=\mD u\by \cdots \mD^pu$ be a $p$th order Jet of $u$ arising along matrix-indexed infinitesimal subsequences 
\[
\de_{D^{[p],\um}u} \, \overset{^*}{\smash{\lharpoonup}}\ \mD^{[p]} u\ \ \   \text{ in }\, \mY\Big(\Om, \smash{\overline{\R}}^{Nn} \by ...\by \smash{\overline{\R}}^{Nn^p}_s\Big),
\]
as $\um \ri \infty$ (Definitions  \ref{definition6}, \ref{Multi-indexed convergence},  \ref{Diffuse Derivatives}, \ref{definition13}). Since the weak* topology on the Young measures is metrisable, we may apply Step 1 to 
\[
\mathfrak{X}\, =\, \mY\Big(\Om, \smash{\overline{\R}}^{Nn} \by ...\by \smash{\overline{\R}}^{Nn^p}_s\Big),\ \ \ f^\um=\,  \de_{D^{[p],\um}u},\ \ \ f\,=\, \mD^{p]}u,
\]
for some metric $\rho$ inducing the weak* topology. Thus, we infer that there exists a single-indexed subsequence $(\um_{\nu})_1^\infty \sub \N^{p\by p}$ such that
\[
 \lim_{\nu\ri\infty} \, \rho\left(  \de_{D^{[p],\um_\nu }u} ,\, \mD^{p]}u\right) \,=\,0,
\]
because by the definition of $\mD^{[p]} u$ we have
\[
\lim_{m^p_p \ri\infty} \left(... \, ...\lim_{m^2_2 \ri\infty}\, \left(\lim_{m^1_2 \ri\infty} \left( \lim_{m^1_1 \ri\infty} \, \rho\left(  \de_{D^{[p],\um}u} ,\, \mD^{p]}u\right)\right)\right)\right)\,=\,0.
\]
Since $u$ is measurable, by invoking Lemma \ref{lemma10} for $\e=1/|\um| $ we obtain a multi-indexed sequence $(u^{\underline{m}})_{\underline{m} \in \N^{p\by p}} \sub C^\infty_0(\Om,\R^N)$. Let $(u^\nu)_{\nu \in \N}$ be the subsequence of it corresponding to $(\um_{\nu})_1^\infty$. Then, by \eqref{3.1} and by recalling that almost uniform implies a.e.\ convergence, we immediately have $u^\nu \larrow u$ a.e.\ on $\Om$ as $\nu \ri \infty$. Further, again by \eqref{3.1} we have
\[
\Big| D^{[p],\um_\nu }u \, -\, D^{[p]}u^\nu\Big|\, \larrow \,0 \ \ \text{ a.e.\ on }\Om,\ \text{ as }\nu\ri\infty.
\]
By Lemma \ref{lemma34} and the above, we obtain
\[
\de_{D^{[p]}u^\nu} \, \overset{^*}{\smash{\lharpoonup}}\ \mD^{[p]} u\ \ \   \text{ in }\, \mY\Big(\Om, \smash{\overline{\R}}^{Nn} \by ...\by \smash{\overline{\R}}^{Nn^p}_s\Big), \ \text{ as }\nu\ri\infty.
\]

\smallskip

\noi  \textbf{Step 3}.  We now \emph{define} for each $\nu \in \N$ the mapping $f^\nu :\Om\sub \R^n \larrow\R^M$ given by
\[
f^\nu(x)\,:=\, F\Big(x,u^\nu(x),D^{[p]}u^\nu(x) \Big) , \ \ \ \ x\in \Om.
\]
In order to conclude the theorem, we seek to show that for any 
\[
\Phi\, \in C^0_c\big( \R^{Nn} \! \by \cdots \by \R^{Nn^p}_s \big), 
\]
we have that
\[
\text{for a.e.\ $x\in \Om$, \ \ $\Phi\Big( D^{[p]}u^\nu(x)\Big) f^\nu(x) \larrow 0$ \ as $\nu \ri\infty$.}
\]
This last statement is a consequence of the Convergence Lemma \ref{lemma16A}. Indeed, since $u$ is a $\mD$-solution on $\Om$, for a.e.\ $x\in \Om$, we have
\[ 
\left| F\big( x,u(x),\underline{\X}\big)\right| \, =\, 0, \ \ \ \ \uX \in \, \supp_*\big( \mD^{[p]} u (x)\big).
\]
We fix such an $x$ as above and we choose a function $\Phi$ as above. Then, we note that the definition of $\mD$-solutions implies that the continuous function
\[
\smash{\overline{\R}}^{Nn} \by\cdots \by \smash{\overline{\R}}^{Nn^p}_s \ \ni\ \uX \lmapsto \left| \Phi(\underline{\X}) \, F\big( x,u(x),\underline{\X}\big)\right| \ \in \ \R \ \ \ \ 
\]
is well-defined on the compactification and vanishes on the support of the probability measure $\mD^{[p]} u (x)$. Hence, we have
\[
\int_{\smash{\overline{\R}}^{Nn} \by\cdots \by \smash{\overline{\R}}^{Nn^p}_s}  \left| \Phi(\underline{\X}) \, F\big( x,u(x),\underline{\X}\big)\right| d\big[\mD^{[p]} u (x)\big](\underline{\X})\, =\, 0, \quad \text{ a.e. }x\in\Om.
\]
Since $u^\nu \larrow u$ a.e.\ on $\Om$ and $\de_{D^{[p]}u^\nu} \, \overset{^*}{\smash{\lharpoonup}}\ \mD^{[p]} u$, by applying Lemma \ref{lemma16A} we obtain that  
\[
\begin{split}
 \Big| \Phi\big( D^{[p]}u^\nu \big) f^\nu\Big|(x) \, &=\,  \left| \Phi \left(D^{[p]}u^\nu (x)\right) \, F\Big( x,u^\nu (x),D^{[p]}u^\nu (x)\Big)\right| 
\\
&= \, \int_{\smash{\overline{\R}}^{Nn} \by\cdots \by \smash{\overline{\R}}^{Nn^p}_s}  \left| \Phi(\underline{\X}) \, F\big( x,u^\nu (x),\underline{\X}\big)\right| d\big[\de_{D^{[p]}u^\nu (x)}\big](\underline{\X}) 
\\
& \!\!\! \larrow  0, 
\end{split}
\]
as $\nu\ri\infty$, for a.e.\ $x\in \Om$. The theorem follows.  \qed

\ms

Now is remains to establish Corollary \ref{corollary12}.

\BPCOR \ref{corollary12}. We continue from the proof of Theorem \ref{theorem11}. First we note that the equivalences among the modes of convergence described in \eqref{3.30}-\eqref{3.31} are quite elementary. By choosing $\Phi \geq \chi_{\mB_R(0)}$ or $\Phi \leq C\chi_{\mB_R(0)}$ we see that \eqref{3.30} is equivalent to \eqref{3.31}(a). Then, \eqref{3.31}(a) is equivalent to \eqref{3.31}(b) by noting that a.e.\ convergence is equivalent (up to the passage to a subsequence) to convergence locally in measure and also by observing the identity
\[
\Big\{ |f^\nu|\, \chi_{\mB_{1/\e}(0)}\big( D^{[p]}u^\nu\big) >\e  \Big\}\, =\, \big\{ |f^\nu|>\e\big\} \cap \Big\{ \big|D^{[p]}u^\nu\big|<{1}/{\e}\Big\} 
\]
which is valid for any $\e>0$. Hence, in order to conclude it suffices to establish \eqref{3.32} when the set $E\sub \Om$ is given by \eqref{3.33}. To this end, we fix $\Om'\Subset \Om$ and we consider the sequence
\[
\psi^\nu(x)\, :=\, \min\big\{|f^\nu(x)|\, \chi_{\Om'\set E}(x),\, 1 \big\}, \ \ \ \ x\in \Om, \ \nu\in \N.
\]
Note now that $(\psi^\nu)_1^\infty$ is bounded and equi-integrable in $L^1(\Om)$ because
\[
\lim_{t\ri\infty} \left(\sup_{\nu\in \N}\, \int_{\Om\cap \{|\psi^\nu|>t\}} |\psi^\nu(x)|\,dx \right)=\, 0.
\]
Moreover, the integral
\[
\int_{\Om'\set E} \int_{ { {\R}}^{Nn} \by\cdots \by  { {\R}}^{Nn^p}_s}  \min \Big\{ \big| F\big( x,u (x),\underline{\X}\big)\big|, 1\Big\} \, d\big[\mD^{[p]}u (x)\big](\underline{\X}) \, dx
\]
exists because for a.e.\ $x\in \Om\set E$, we have that $\big[\mD^k u(x)\big](\{\infty\})=0$ for all $k=1,...,p$. In addition, by utilising that $u^\nu\larrow u$ a.e.\ on $\Om\set E$ and Lemma \ref{lemma35}, we have
\[
\de_{(u^\nu, D^{[p]}u^\nu)} \, \overset{^*}{\smash{\lharpoonup}}\ \de_u \!\by \mD^{[p]} u\ \ \   \text{ in }\, \mY\Big(\Om\set E,\,  {\R}^{N}\! \by {\R}^{Nn} \by ...\by  {\R}^{Nn^p}_s\Big),
\]
as $\nu\ri\infty$. Hence, by invoking Corollary 3.36 on p.\ 207 of \cite{FG}, we have
\begin{align}  \label{+}
\ \ \ \ \ \ \lim_{\nu\ri\infty}   \int_{\Om'\set E}& \big|\psi^\nu  (x)  \big|\,dx \nonumber \\
 =\, &
\lim_{\nu\ri\infty} \int_{\Om'\set E} \min \Big\{ \Big| F\Big( x,u^\nu(x),D^{[p]}u^\nu(x)\Big)\Big|, 1\Big\} \,dx\\
=\  & 
\int_{\Om'\set E} \int_{ { {\R}}^{Nn} \by\cdots \by  { {\R}}^{Nn^p}_s}  \min \Big\{ \big| F\big( x,u (x),\underline{\X}\big)\big|, 1\Big\} \, d\big[\mD^{[p]}u (x)\big](\underline{\X}) \, dx.  \nonumber
\end{align}
Since $u$ is a $\mD$-solution, we have
\[
\sup_{ \uX \in\, \supp_* (\mD^{[p]}u (x) ) } \big| F\big( x,u (x),\underline{\X}\big)\big|\, =\,0\ \ \ \ \text{ a.e. on }\Om.
\]
Moreover, for any $\e\in (0,1)$ we have the inequality
\[
\begin{split}
\e\, \Big|  (\Om'\set E)  \cap  \big\{  |f^\nu|  >\e\big\}\Big| \,
 &\leq  \,  \int_{\Om'\set E} \min\big\{ |f^\nu(x)|,1\big\}\, dx
 \\
 &=\,  \int_{\Om} \big|\psi^\nu  (x)  \big|\,dx
  \end{split}
\]
Further, we note that for a.e.\ $x\in \Om\set E$, $\mD^{[p]}u (x)$ is a probability measure on ${\R}^{Nn} \! \by \cdots \by  {\R}^{Nn^p}_s $ (not just on the compactification). By recalling the definition of the function $\psi^\nu$, \eqref{+} and the above observations give
\[
\begin{split}
\e\, \underset{\nu\ri\infty}{\lim\sup} \, &\Big|  (\Om'\set E)  \cap  \big\{  |f^\nu|  >\e\big\}\Big| \\
\leq   &\, 
\int_{\Om'\set E} \int_{ { {\R}}^{Nn} \by\cdots \by  { {\R}}^{Nn^p}_s}  \min \Big\{ \big| F\big( x,u (x),\underline{\X}\big)\big|, 1\Big\} \, d\big[\mD^{[p]}u (x)\big](\underline{\X}) \, dx\\
 \leq &\,  
\int_{\Om'\set E} \int_{ { {\R}}^{Nn} \by\cdots \by  { {\R}}^{Nn^p}_s} \big| F\big( x,u (x),\underline{\X}\big)\big| \, d\big[\mD^{[p]}u (x)\big](\underline{\X}) \, dx\\
 \leq & \, 
\int_{\Om'\set E} \left(\sup_{ \uX \in\, \supp \big(\mD^{[p]}u (x)\big) \cap \big( \R^{Nn} \by\cdots \by  \R^{Nn^p}_s \big) } \big| F\big( x,u (x),\underline{\X}\big)\big| \right) dx\\
=&\, 
\int_{\Om'\set E} \left(\sup_{ \uX \in\, \supp_* (\mD^{[p]}u (x) ) } \big| F\big( x,u (x),\underline{\X}\big)\big| \right) dx\\
 =&\,  \ 0.
\end{split}
\]
Conclusively, we have obtained that $f^\nu \larrow 0$ locally in measure on $\Om\set E$ and hence up to a subsequence $f^\nu \larrow 0$ a.e.\ on $\Om\set E$ as $\nu\ri\infty$. The corollary has been established.         \qed

\ms
\ms

\noi \textbf{Acknowledgement.} The author would like to thank N. Barron, J. Manfredi and P. Pedregal for their comments and encouragement. 

\ms

\bibliographystyle{amsplain}

\end{document}